\documentclass[a4paper,10pt]{article}
\usepackage{amssymb,amsthm,amsmath}    
\usepackage[english]{babel}                         
\usepackage[ansinew]{inputenc}         
\usepackage[dvips]{graphicx}           
\theoremstyle{definition}
\newtheorem{exam}{Example}

\title{{\small For Dagstuhl Seminar\\ Combinatorial and Algorithmic Aspects of Sequence Processing 11081\\Organizers: Maxime Crochemore, Lila Kari, Mehryar Mohri and Dirk Nowotka\\Date: 20.02.2011-25.02.2011} \\ Observations and Problems on $k$-abelian avoidability \footnote{Supported by the Academy of Finland under the grant 121419 and by the Väisälä Foundation.}}
\author{Mari Huova and Juhani Karhumäki} 
\date{{\small Department of Mathematics and TUCS \\ University of Turku \\  20014 Turku, FINLAND \\ email:  $\left\{ \textrm{mari.huova, karhumak} \right\}$@utu.fi} }

\begin{document}

\maketitle

Theory of avoidability is among the oldest and most studied topic in Combinatorics on Words. The first result in this area, or in fact in the whole field, were obtained by Norwegian Axel Thue as early as at the beginning of 20th century \cite{Thue1, Thue2}. He showed, among other things, the existence of an infinite binary word, which does not contain three consecutive factors of a word, that is a cube. Similarly, he showed that squares can be avoided in infinite ternary words.

Since late 1960's commutative variants of the above problems were studied. Apparently, first nontrivial results were obtained by Evdokimov \cite{Evdok} who showed that commutative squares could be avoided in infinite words over a 25-letter alphabet. The size of the alphabet was reduced to five by Pleasant \cite{Pleas}, until the optimal value four was found by Keränen \cite{Ker}, solving one celebrated problems of the topic. The optimal value for the size of the alphabet where abelian cubes were avoidable was proved earlier by Dekking \cite{Dek}, the value being three.

Interesting in all these results is that the required words are obtained by iterating a morphism.

We introduce in this note new variants of the problems by defining repetitions via new equivalence relations which lie properly in between equality and commutative equality, that is abelian equality.

Let $k \geq 1$ be a natural number. We say that words $u$ and $v$ in $\Sigma^{+}$ are \textit{$k$-abelian equivalent}, in symbols $u \equiv_{a,k} v$, if

\begin{enumerate}
	\item $\textrm{pref}_{k-1} \left( u \right) = \textrm{pref}_{k-1} \left( v \right)$ and $\textrm{suf}_{k-1} \left( u \right) = \textrm{suf}_{k-1} \left( v \right)$, and
	\item for all $w \in \Sigma^{k}$, the number of occurrences of $w$ in $u$ and $v$ coincide, i.e. $\#\left( w,u \right) = \#\left( w,v \right)$.
\end{enumerate}

\noindent Here $\textrm{pref}_{k-1}$ (resp. $\textrm{suf}_{k-1}$) is used to denote the prefixies (resp. suffixies) of length $k-1$ of words.

It is straightforward to see that $\equiv_{a,k}$ is an equivalence relation and, moreover,
$$ u=v \Rightarrow u \equiv_{a,k} v \Rightarrow u \equiv_{a} v, $$

\noindent where $\equiv_{a}$ denotes the abelian equivalence, and that
$$ u=v \Leftrightarrow u \equiv_{a,k} v \ \ \forall\; k \geq 1. $$

Now, notions like \textit{$k$-abelian repetitions} are naturally defined. For instance, $w = uv$ is a \textit{$k$-abelian square} if and only if $u \equiv_{a,k} v$.

Natural variants of the above Thue's problems ask what are the smallest alphabets where $k$-abelian square and cubes can be avoided. A goal of this note is to point out that these problems are not trivial even in the case $k = 2$. Before going into that we make a few preliminary simple observations. 

First, in the binary case 2- and 3-abelian words are fairly easy to characterize.

\begin{exam}
In a binary alphabet $\Sigma = \left\{a,b\right\}$ the characterization of equivalence classes of 2-abelian words, via their representatives, can be given in the form:
$$aa^{k}b^{l}(ab)^{m}a^{n} \  \textrm{ or } \  bb^{k}a^{l}(ba)^{m}b^{n} ,$$
where $k,l,m \geq 0 $ and $n \in \left\{ 0,1 \right\}$.
And in the same alphabet the characterization of equivalence classes of 3-abelian words can be given in the following form containing eight possible combinations:
\begin{displaymath}
\left. \begin{array}{lll}
\left. \begin{array}{l}
aaa^{k}b^{l}(aabb)^{m}  \\
bbb^{k}a^{l}(aabb)^{m} \\ 
abb^{k}a^{l}(aabb)^{m}  \\ 
baa^{k}b^{l}(aabb)^{m}  \\
\end{array} \textrm{or} \right\} \ast & 
\textrm{connected with} &

\ast \left\{ \begin{array}{l}
(aab)^{g}(ab)^{h}b^{i}a^{j} \\
(abb)^{g}(ab)^{h}b^{i}a^{j}, \\
\end{array} \textrm{or} \right.\\
\end{array} \right.
\end{displaymath}
where $k,l,m,g,h \geq 0 $ , $i \in \left\{ 0,1 \right\}$ and $j \in \left\{ 0,\ldots,2-i \right\}$.
Here the characterizations are not unique in few cases.
\end{exam}

The above allows to estimate the sizes of the corresponding equivalence classes. They are of order $\Theta \left( n^{2} \right)$ and $\Theta \left( n^{4} \right)$, see \cite{Hkss}. Recently, A. Saarela \cite{Saar} showed that in general the number of $k$-abelian equivalence classes of words of length $n$ is polynomial in $n$ but the degree of the polynomial grows exponentially in $k$ (in a fixed but arbitrary alphabet).

Our next example shows that the ordinary method of iterating a morphism might not give answers to our problems.

\begin{exam} In all of the following cases where repetition free infinite word is obtained by iterating a morphism, a 2-abelian cube is found fairly early from the beginning. For overlap- and cube-free words see \cite{Aj}.

\begin{itemize}
	\item Infinite overlap-free Thue-Morse word (morphism: $0\rightarrow01$, $1\rightarrow10$):
$01\overbrace{101001}\overbrace{100101}\overbrace{101001}011...$
  \item Cube-free infinite word (morphism: $0\rightarrow001$, $1\rightarrow011$):
$001001\overbrace{011001}\overbrace{001011}\overbrace{001011}011...$
	\item Morphism $0\rightarrow001011$, $1\rightarrow001101$, $2\rightarrow011001$ maps ternary cube-free words to binary cube-free words, see \cite{Bran}, but $001011 \equiv_{a,2} 001101 \equiv_{a,2} 011001$, thus images of all words mapped with this morphism contains 2-abelian cubes.
	\item A binary overlap-free word $w$ can also be gained in form $w= c_{0}c_{1}c_{2} \ldots$, where $c_{n}$ means the number of zeros (mod 2) in the binary expansion of $n$. Again, a 2-abelian cube of length 6 begins as early as from the fifth letter: $w= 0010\overbrace{011010}\overbrace{010110}\overbrace{011010}011...$
\end{itemize}
\end{exam}

In order to go into our problems we recall the following Table \ref{tab1} which summarizes the results we mentioned at the beginning and at the same time tells the limits of our problems:

\begin{table}[h!]
\begin{tabular}{|c|c|c|c||c|c|c|c|}
\hline
\multicolumn{4}{|c||}{Avoidability of squares} & \multicolumn{4}{|c|}{Avoidability of cubes}\\
\hline
 & \multicolumn{3}{|c||}{type of rep.} &  & \multicolumn{3}{|c|}{type of rep.}\\
size of the alph. &  $=$ & $\equiv_{a,2}$ & $\equiv_{a}$ & size of the alph. &  $=$ & $\equiv_{a,2}$ & $\equiv_{a}$ \\
\hline
2 & $-$ &  $-$ &  $-$ & 2 & $+$ &  $?$ &  $-$\\
\hline
3 & $+$ & ? &  $-$  & 3 & $+$ & $+$ &  $+$\\
\hline
4 & $+$ &  $+$ &  $+$ & \multicolumn{4}{|c|}{}\\
\hline
\end{tabular}
\caption{Avoidability of different types of repetitions in infinite words.}\label{tab1}
\end{table}

We were able to settle the first one of the above question marks by computer checking.

\begin{exam}
The longest ternary word which is 2-abelian square-free has length 537, which shows that there does not exist an infinite 2-abelian square-free word over any ternary alphabet. This longest word, given below, is unique up to the permutations of the alphabet, $\Sigma = \left\{ a,b,c \right\}$. 

\begin{displaymath}
\left. \begin{array}{l}
abcbabcacbacabacbabcbacabcbabcabacabcacbacabacbabcbacbcacbabcacbcabcba\\
bcabacbabcbacbcacbacabacbabcbacabcbabcabacabcacbacabacbabcbacbcacbacab\\
acbcabacabcacbcabcbacbcacbacabacbabcbacbcacbabcacbcabcbabcabacbabcbacb\\
cacbacabacbabcbacabcbabcabacabcacbacabacbabcbacbcacbacabacbcabacabcacb\\
cabcbabcabacabcacbacabacbabcbacabcbabcabacabcacbcabcbabcabacbabcbacbca\\
cbabcacbcabcbabcabacabcacbcabcbacbcacbacabacbcabacabcacbcabcbabcabacab\\
cacbacabacbabcbacabcbabcabacabcacbcabcbabcabacbabcbacbcacbabcacbcabcba\\
bcabacabcacbacabacbabcbacabcbabcabacabcacbabcba\\

\end{array} \right.
\end{displaymath}

\end{exam}

\begin{exam}
We can also construct the whole tree containing each ternary 2-abelian square-free word and analyze the sizes of the sets containing words from length 1 to 537, respectively. There exist 404 286 words of length 105 and for other lengths the number of ternary 2-abelian square-free words is less. The number of words grows monotonically from 3 up to 403 344 when considering lengths from 1 to 103. After the length 105 there appears more oscillation between the numbers of words of different lengths. The sizes of the sets containing ternary 2-abelian square-free words are shown in Figure \ref{pic1} with respect to the lengths of words.

\begin{figure}[h]
\includegraphics[width=12cm]{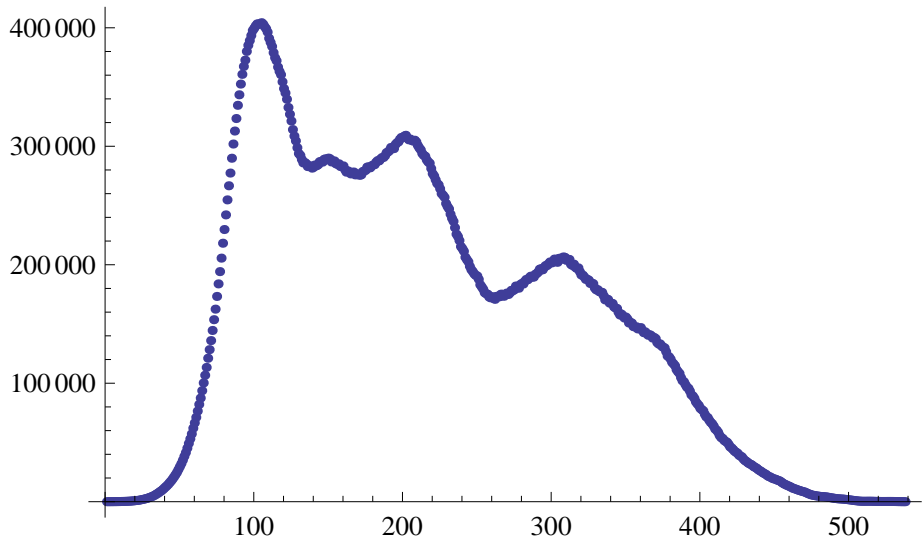}
\caption{The number of 2-abelian square-free words with respect to their lengths.} \label{pic1}
\end{figure}
\end{exam}
To solve the other question mark we also did some computer checking - and obtained evidence that the answer is likely to be different compared to the first one.

\begin{exam}
With a computer we were able to construct a binary word of more than 100 000 letters that still avoids 2-abelian cubes. This shows that there exist, at least, very long binary 2-abelian cube-free words.
\end{exam}

\begin{exam}
Similarly, we can examine the number of binary 2-abelian cube-free words of given length as in the previous case concerning ternary 2-abelian square-free words. The numbers of the words with lengths from 1 to 60 grow approximately with a factor 1,3 at each increment of the length, see Figure \ref{pic2}. So that the number of binary 2-abelian cube-free words of length 60 is already 478 456 030. And already, with length 12 there exist more binary 2-abelian cube-free words (254) than ternary 2-abelian square-free words (240).

\begin{figure}[h!]
\includegraphics[width=12cm]{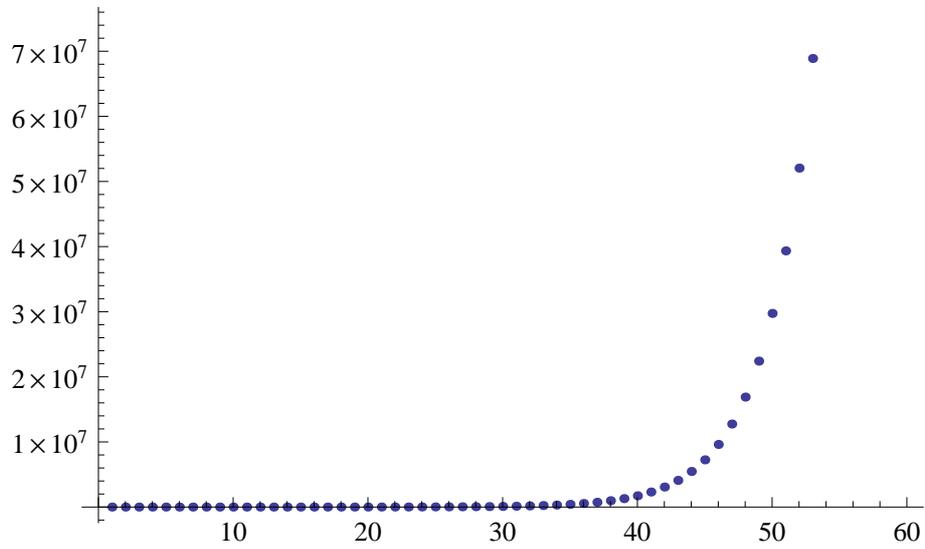}
\caption{The number of 2-abelian cube-free words with respect to their lengths for small values of length.} \label{pic2}
\end{figure}

We also chose some binary 2-abelian cube-free prefixies and counted the numbers of binary 2-abelian cube-free words having these fixed preixies. In this way we can check how many suitable extensions the chosen 2-abelian cube-free word has. As a result we found examples of binary 2-abelian cube-free words with a property that the number of their extensions grows again approximately with a factor 1,3 when increasing the length of extensions by one. In Figure \ref{pic3} this is done fore a fixed prefix of length 2000. 

\begin{figure}[h!]
\includegraphics[width=12cm]{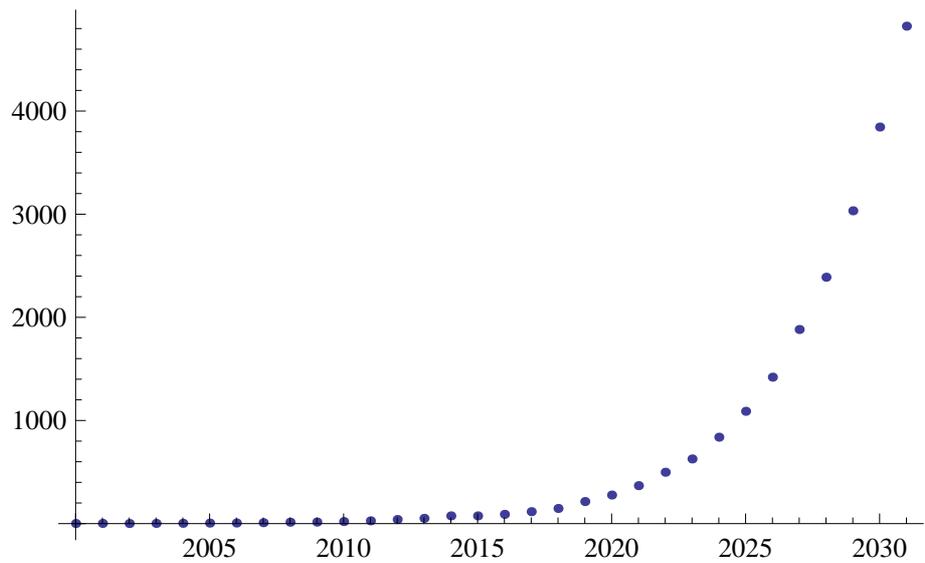}
\caption{The number of 2-abelian cube-free words with respect to their lengths from 2000 to 2031 and with a common prefix of length 2000.} \label{pic3}
\end{figure}
\end{exam}

These examples support the conjecture that there would exist an infinite binary word that avoids 2-abelian cubes. As a conclusion, our two considered problems would behave differently: one like words and the other like abelian words.

\renewcommand{\refname}{References}


\begin{thebibliography}{99}
\bibitem[AJ]{Aj} J.-P. Allouche, J. Shallit: \textsl{Automatic Sequences: Theory, Applications, Generalizations}. Cambridge University Press, Cambridge, 2003.

\bibitem[Br]{Bran} F.-J. Brandenburg: \textsl{Uniformly growing k-th power-free homomorphisms}. Theoret. Comput. Sci. 23, 69-82 (1983).

\bibitem[De]{Dek} F. M. Dekking: \textsl{Strongly non-repetitive sequences and progression-free sets}. J. Combin. Theory Ser. A 27(2), 181-185 (1979).

\bibitem[Ev]{Evdok} A. A. Evdokimov: \textsl{Strongly asymmetric sequences generated by a finite number of symbols}. Dokl. Akad. Nauk SSSR 179, 1268-1271 (1968); English translation in Soviet Math. Dokl. 9, 536-539 (1968).

\bibitem[HKSS]{Hkss} M. Huova, J. Karhumäki, A. Saarela, K. Saari: \textsl{Local squares, periodicity and finite automata}. LNCS Festschrift for Hermann Maurer, Springer, (to appear).

\bibitem[Ke]{Ker} V. Keränen: \textsl{Abelian squares are avoidable on 4 letters}. In: W. Kuich (ed.) ICALP 1992. LNCS, vol. 623, 41-52. Springer, Heidelberg, 1992.

\bibitem[Pl]{Pleas} P. A. B. Pleasant: \textsl{Non-repetitive sequences}. Proc. Cambridge Philos. Soc. 68, 267-274 (1970).

\bibitem[Sa]{Saar} A. Saarela: Private communication.

\bibitem[Th1]{Thue1} A. Thue: \textsl{Über unendliche Zeichenreihen}. Norske vid. Selsk. Skr. Mat. Nat. Kl. 7, 1-22 (1906).

\bibitem[Th2]{Thue2} A. Thue: \textsl{Über die gegenseitige Lage gleicher Teile gewisser Zeichenreihen}. Norske vid. Selsk. Skr. Mat. Nat. Kl. 1, 1-67 (1912).

\end{thebibliography}
\end{document}